\input amstex
\documentstyle{amsppt}
\pagewidth{5.4in}
\pageheight{7.6in}
\magnification=1200
\TagsOnRight
\NoRunningHeads
\topmatter
\title
\bf A lower bound for the scalar curvature of\\
the standard solution of the Ricci flow  
\endtitle
\author
Shu-Yu Hsu
\endauthor
\affil
Department of Mathematics\\
National Chung Cheng University\\
168 University Road, Min-Hsiung\\
Chia-Yi 621, Taiwan, R.O.C.\\
e-mail:syhsu\@math.ccu.edu.tw
\endaffil
\date
Dec 15, 2006
\enddate
\address
e-mail address:syhsu\@math.ccu.edu.tw
\endaddress
\abstract
In this paper we will give a rigorous proof of the lower bound for the 
scalar curvature of the standard solution of the Ricci flow conjectured by 
G.~Perelman. We will prove that the scalar curvature $R$ of the standard 
solution satisfies $R(x,t)\ge C_0/(1-t)\quad\forall x\in\Bbb{R}^3,0\le t<1$, 
for some constant $C_0>0$. 
\endabstract
\endtopmatter
\NoBlackBoxes
\define \1{\partial}
\define \2{\overline}
\define \3{\varepsilon}
\document

Recently there is a lot of study of Ricci flow on manifolds by R.~Hamilton
\cite{H1-6}, S.Y.~Hsu \cite{Hs1-6}, P.~Lu and G.~Tian \cite{LT},
G.~Perelman \cite{P1}, \cite{P2}, W.X.~Shi \cite{S1}, \cite{S2}, L.F.~Wu
\cite{W1}, \cite{W2}, and others. In \cite{H1} R.~Hamilton studied the 
Ricci flow on compact manifolds with strictly positive Ricci curvature. 
He proved that if the metric $g(x,t)$ of the manifold evolves 
by the Ricci flow,
$$
\frac{\1 }{\1 t}g_{ij}=-2R_{ij}\tag 0.1
$$
with $g_{ij}(x,0)=g_{ij}(x)$, then the evolving metric will converge modulo 
scaling to a metric of constant positive curvature. This result was later
extended to compact four dimensional manifold with positive curvature operator
by R.~Hamilton \cite{H2} and to non-compact complete manifolds by W.X.~Shi 
\cite{S1}, \cite{S2}. Behaviour and properties of Ricci flow on $\Bbb{R}^2$ 
are also studied by P.~Daskalopoulos and M.A.~Del Pino \cite{DP}, S.Y.~Hsu 
\cite{Hs1-4} and L.F.~Wu \cite{W1}, \cite{W2}. We refer the reader to the 
survey paper of R.~Hamilton \cite{H5} for previous results on Ricci flow,  
the lecture notes \cite{Ch} of B.~Chow and the book \cite{CK} of B.~Chow 
and D.~Knopf for the most recent results on Ricci flow.
 
In \cite{P1}, \cite{P2}, G.~Perelman introduced the concepts of
$L$-geodesic, $L$-length, and reduced volume, to study the Ricci flow on
manifolds with singularities. In \cite{P2} G.~Perelman studied the Ricci 
flow with surgery on manifolds. Essential to this program is the 
construction of a standard solution for the Ricci flow on $\Bbb{R}^3$ with 
certain properties. This standard solution is then used to replace the 
solution of the Ricci flow near singularity points during surgery. In 
\cite{P2} G.~Perelman conjectured that the scalar curvature $R(x,t)$ of the 
standard solution satisfies
$$
R(x,t)\ge\frac{C_0}{1-t}\quad\forall x\in\Bbb{R}^3,0\le t<1\tag 0.2
$$
for some constant $C_0>0$. However there is no detailed proof of this result 
in the notes on Perelman's papers \cite{KL} by B.~Kleiner and J.~Lott. The 
proof of this result in \cite{CZ} (Proposition 7.3.3) by H.D.~Cao and 
X.P.~Zhu is questionable. The proof of this result in 
\cite{MT} (Proposition 12.31) by J.~Morgan and G.~Tian is incomplete and 
has gaps. In this paper we will give a simple  correct proof of (0.2). 

We first start with a definition. Let $\kappa>0$. A Ricci flow $(M,g)$ is 
said to be $\kappa$-noncollapsing at the point $(x_0,t_0)$ on the scale 
$r>0$ \cite{P1} if 
$$
\text{Vol}_{g(t_0)}(B_{g(t_0)}(x_0,r))\ge\kappa r^n
$$
holds whenever
$$
|\text{Rm}|(x,t)\le r^{-2}\quad\forall d_{g(t)}(x_0,x)<r, t_0-r^2\le 
t\le t_0
$$
holds where $B_{g(t_0)}(x_0,r)$ is the geodesic ball of radius $r$ in $M$ 
around the point $x_0$ with respect to the metric $g(t_0)$. A Ricci flow 
$(M,g)$ is said to be a $\kappa$-solution if it is a solution of the Ricci 
flow in $M\times (-\infty,0)$ such that for each $t<0$ the metric $g(t)$ is
not a flat metric, $(M,g(t))$ is a complete manifold of nonnegative 
curvature, $M$ has uniformly bounded curvature in $(-\infty,0)$ and 
$(M,g)$ is $\kappa$-noncollapsing at all points of $M\times 
(-\infty,0)$. 

For any $0\le t<1$, let $h(t)$ be the standard metric on $S^2$ whose scalar 
curvature on $S^2$ is $1/(1-t)$. Let $g_0$ be a fixed rotationally symmetric 
complete smooth metric with positive curvature operator on $\Bbb{R}^3$
such that $(\Bbb{R}^3\setminus\2{B(0,2)},g_0)$ is isometric to the half
infinite cylinder $(S^2\times\Bbb{R}^+,h(1)\times ds^2)$ (cf \cite{MT} 
definition 12.1). By Lemma 12.2 of \cite{MT} such $g_0$ exists. We say that
a Ricci flow $(\Bbb{R}^3,g(t))$, $0\le t<1$, is a standard solution if 
$g(0)=g_0$ and the curvature $Rm$ is locally bounded in time $t\in [0,1)$. 
By the results of \cite{P2}, \cite{KL}, \cite{S1} and \cite{S2}, there exists 
a unique standard solution $(\Bbb{R}^3,g(t))$ of the Ricci flow on $(0,1)$
with $g(0)=g_0$ which has positive curvature operator $Rm(t)$ for each 
$t\in [0,1)$. 

$$
\text{Section 1}
$$

We first recall some results of \cite{P1}, \cite{P2}, 
\cite{KL}, and \cite{MT}.

\proclaim{\bf Theorem 1.1}(cf. Theorem 12.1 of \cite{P1} and Theorem 51.3 of
\cite{KL}) For any $\3>0$, $\kappa>0$ and $\sigma>0$, there exists a constant 
$r_0>0$ such that the following holds. Let $T>1/2$. Suppose $g(t)$, 
$0\le t<T$, is a solution to the Ricci flow on a three-manifold $M$ such 
that for each $0\le t<T$, $(M,g(t))$ is a complete manifold of bounded 
sectional curvature. Suppose the Ricci flow also has 
non-negative curvature and is $\kappa$-noncollapsing on scales less than 
$\sigma$. Then for any $(x_0,t_0)\in M\times [1/2,T)$ satisfying $Q=R
(x_0,t_0)\ge r_0^{-2}$, the solution in $\{(x,t):\text{dist}_{t_0}^2(x_0,x)
<(\3Q)^{-1}, t_0-(\3Q)^{-1}\le t\le t_0\}$ is after scaling by the factor 
$Q$, $\3$-close to the corresponding subset of a $\kappa$-solution.   
\endproclaim

\proclaim{\bf Lemma 1.2}(Lemma 58.4 of \cite{KL} and Section 1.5 of \cite{P2})
Let $\kappa>0$. Then there exists a constant $C_1>0$ such that for any 
$\kappa$-solution $(M,g)$ the scalar curvature satisfies,
$$
|\nabla R|\le C_1R^{\frac{3}{2}},|R_t|\le C_1R^2.\tag 1.1
$$
for any $(x,t)\in M\times (-\infty,0)$.
\endproclaim

We will now assume that $(\Bbb{R}^3,g(t))$ is the standard solution 
of the Ricci flow for the rest of this paper. For any $t\in [0,1)$, 
$x_0\in\Bbb{R}^3$ and $r>0$, let $B_t(x_0,r)=\{x\in\Bbb{R}^3:d_t(x,x_0)<r\}$ 
and $B_t(r)=B_t(0,r)$ where $d_t(x_0,x)$ is the distance between $x_0$ and $x$ 
with respect to the metric $g(t)$.

\proclaim{\bf Lemma 1.3}(Section 2 of \cite{P2}, Lemma 59.3 of 
\cite{KL}, and Theorem 12.5 of \cite{MT})
For each $0\le t<1$, $g(t)$ is a rotationally symmetric complete metric.
Moreover there exists $\sigma>0$ and $\kappa>0$ such that $(\Bbb{R}^3,g(t))$, 
$0\le t<1$, is $\kappa$-noncollapsed on scales less than $\sigma$ on 
$\Bbb{R}^3\times (0,1)$.
\endproclaim

\proclaim{\bf Lemma 1.4}(cf. P.321--322 of \cite{MT})
Let $(\Bbb{R}^3,g(t))$, $0\le t<1$, be the standard solution of Ricci flow. 
Then 
$$
\liminf_{t\to 1}R(x,t)=\infty\quad\forall x\in\Bbb{R}^3.
$$
\endproclaim

By Theorem 1.1, Lemma 1.2 and Lemma 1.3 we have
\proclaim{\bf Lemma 1.5}
Let $(\Bbb{R}^3,g(t))$, $0\le t<1$, be the standard solution of Ricci flow. 
Then there exist constants $r_0>0$ and $C_1>0$ such that for any $(x_0,t_0)\in
\Bbb{R}^3\times [1/2,1)$ satisfying $(x_0,t_0)\ge r_0^{-2}$, the scalar 
curvature $R(x_0,t_0)$ of the standard solution satisfies (1.1).
\endproclaim

We are now ready to state and prove the main theorem of the paper.

\proclaim{\bf Theorem 1.6}
Let $(\Bbb{R}^3,g(t))$, $0\le t<1$, be the standard solution of Ricci flow. 
Then there exists a constant $C_0>0$ such that (0.2) holds. 
\endproclaim
\demo{Proof}
Let $r_0>0$ and $C_1>0$ be as given by Lemma 1.5. We choose $t_0\in 
[\frac{1}{2},1)$ such that
$$
r_0^2>2C_1(1-t_0)\quad\Rightarrow\quad r_0^{-2}<\frac{1}{2C_1(1-t_0)}.
$$
We claim that 
$$
R(x,t)\ge\frac{1}{2C_1(1-t)}\quad\forall t_0<t<1.\tag 1.2
$$
Suppose the claim is false. Then there exists $x_1\in\Bbb{R}^3$, $t_0<t_1<1$,
such that 
$$
R(x_1,t_1)<\frac{1}{2C_1(1-t_1)}.
$$
By increasing $t_1$ if necessary we may assume without loss of generality 
that
$$
r_0^{-2}<R(x_1,t_1)<\frac{1}{2C_1(1-t_1)}.\tag 1.3
$$
We now divide the proof of the claim into two cases.

\noindent$\underline{\text{\bf Case 1}}:$ $R(x_1,t)>r_0^{-2}\quad\forall
t_1<t<1$.

Then by (1.3) and Lemma 1.5,
$$\align
|R_t(x_1,t)|\le&C_1R(x_1,t)^2\qquad\qquad\quad\forall t_1<t<1\\
\Rightarrow\quad\frac{1}{R(x_1,t_1)}-\frac{1}{R(x_1,t)}\le&C_1(t-t_1)
\qquad\qquad\quad\,\,\,\forall t_1<t<1\\
\Rightarrow\qquad\qquad\qquad\,\,\frac{1}{R(x_1,t)}
\ge&\frac{1}{R(x_1,t_1)}-C_1(t-t_1)\quad\forall t_1<t<1\\
\ge&2C_1(1-t_1)-C_1(1-t_1)\quad\forall t_1<t<1\\
=&C_1(1-t_1)\qquad\qquad\qquad\quad\forall t_1<t<1\\
\Rightarrow\qquad\qquad\qquad\,\,\, R(x_1,t)
\le&\frac{1}{C_1(1-t_1)}\qquad\qquad\qquad\quad\forall t_1<t<1.\tag 1.4
\endalign 
$$ 
\noindent$\underline{\text{\bf Case 2}}$: There exists $t_2\in (t_1,1)$ such 
that $R(x_1,t_2)=r_0^{-2}$.

For any $t_2<t<1$, either
$$
R(x_1,t)\le r_0^{-2}\tag 1.5
$$
or
$$
R(x_1,t)>r_0^{-2}\tag 1.6
$$
holds. Suppose first (1.6) holds. Then there exists $t_3\in [t_2,t)$ such that
$$
R(x_1,t_3)=r_0^{-2}\quad\text{ and }\quad R(x_1,s)>r_0^{-2}\quad\forall
t_3<s\le t.\tag 1.7
$$
Then by (1.7) and Lemma 1.5,
$$\align
|R_t(x_1,s)|\le&C_1R(x_1,s)^2\quad\forall t_3<s\le t\\
\Rightarrow\quad\frac{1}{R(x_1,t_3)}-\frac{1}{R(x_1,t)}\le&C_1(t-t_3)\\
\Rightarrow\qquad\qquad\qquad\,\,\frac{1}{R(x_1,t)}
\ge&r_0^2-C_1(t-t_3)\\
\ge&r_0^2-C_1(1-t_0)\\
\Rightarrow\qquad\qquad\qquad\,\,\, R(x_1,t)
\le&\frac{1}{r_0^2-C_1(1-t_0)}.\tag 1.8
\endalign 
$$ 
By (1.5) and (1.8),
$$
R(x_1,t)\le(r_0^2-C_1(1-t_0))^{-1}\quad\forall t_2\le t<1.
\tag 1.9
$$
By (1.4) of case 1 and (1.9) of case 2, we get
$$
\limsup_{t\to 1}R(x_1,t)\le\max([C_1(1-t_1)]^{-1},
(r_0^2-C_1(1-t_0))^{-1})<\infty.
$$
This contradicts Lemma 1.4. Hence the claim (1.2) must hold.

We next observe that by the definition of the initial metric $g_0$ for the 
standard solution, 
$$
C_2=inf_{\Bbb{R}^3}R(g_0)>0.
$$
Then by an argument similar to the proof of Lemma 4.13 of \cite{S2},
$$ 
R(x,t)\ge C_2\quad\forall x\in\Bbb{R}^3,0\le t<1.\tag 1.10
$$
Let $C_0=\min (C_2(1-t_0),1/(2C_1))$. Then by (1.2) and (1.10), $R(x,t)$ 
satisfies (0.2) and the theorem follows.
\enddemo

\enddocument